\documentclass[12pt]{article}

\usepackage{amssymb}
\usepackage{amsmath,amsthm}
\usepackage[dvips]{graphicx}
\usepackage{wrapfig}

\newtheorem{thm}{Theorem}[section]

\newtheorem{lem}[thm]{Lemma}

\theoremstyle{definition}
\newtheorem{defn}[thm]{Definition}
\newtheorem{eg}[thm]{Example}
\theoremstyle{remark}

\title{The warping degree of a knot diagram}

\author{Ayaka Shimizu \\ \footnotesize{Department of Mathematics, Osaka City University} \\[-3pt]
\footnotesize{Sugimoto, Sumiyoshi-ku Osaka 558-8585, Japan}\\
\footnotesize{Email: ayakaberrycandy@y2.dion.ne.jp}}

\begin{document}

\maketitle

\begin{abstract}
For an oriented knot diagram $D$, the warping degree $d(D)$ is the smallest number of crossing changes which are needed to obtain 
the monotone diagram from $D$ in the usual way. We show that $d(D)+d(-D)+1$ is less than or equal to the crossing number of $D$. 
Moreover the equality holds if and only if $D$ is an alternating diagram. 
For a knot $K$, we also estimate the minimum of $d(D)+d(-D)$ for all diagrams $D$ of $K$ with $c(D)=c(K)$. 
\end{abstract}

\section{Introduction}

\begin{figure}[h]
\begin{center}
\includegraphics[width=5em]{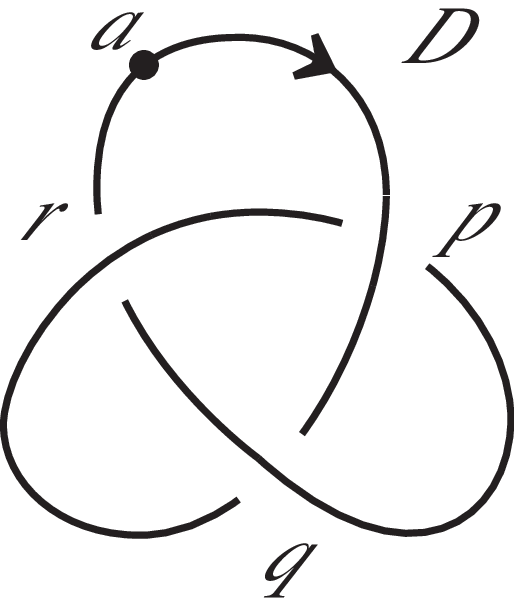}
\caption{}
\label{fig1}
\end{center}
\end{figure}

\medskip
\hspace{3ex} Several concepts of this section are due to Kawauchi \cite{kawauchi}. 
Let $D$ be an oriented knot diagram. 
Let $a$ be a point on $D$ which is not any crossing point. We call it a \textit{base point} of $D$. 
We denote by $D_a$ the pair of $D$ and $a$. 
A crossing point of $D_a$ is a \textit{warping crossing point} if we meet the point first at the under-crossing 
when we go along the oriented diagram $D$ by starting from $a$.
For example, in the diagram in Figure \ref{fig1}, $q$ is a warping crossing point of $D_a$, 
and $p$, $r$ are non-warping crossing points of $D_a$. 
In relation to the warping crossing point, we define the warping degree.

\medskip

\phantom{x}
\begin{defn}
The \textit{warping degree} of $D_a$, denoted by $d(D_a)$, is the number of warping crossing points of $D_a$. 
The \textit{warping degree} of $D$, denoted by $d(D)$, is the minimal warping degree for all base points of $D$. 
\end{defn}
\phantom{x}

Here is an example of $d(D)$. 

\phantom{x}
\begin{eg}
$$
d \left( \begin{minipage}{20pt}
              \includegraphics[width=20pt]{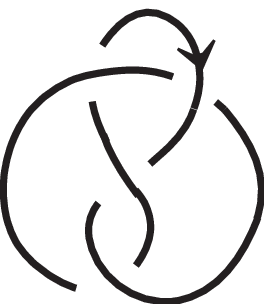}
         \end{minipage}
\right) =1,  \
d \left( \begin{minipage}{20pt}
          \includegraphics[width=20pt]{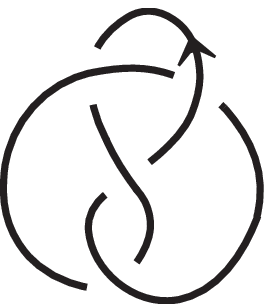}
            \end{minipage}
\right) =2.
$$

\end{eg}
\phantom{x}

\begin{figure}[h]
\begin{center}
\includegraphics[width=5em]{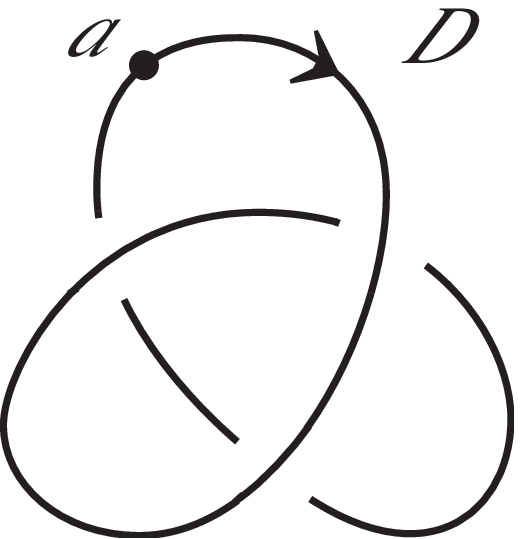}
\caption{}
\label{fig2}
\end{center}
\end{figure}

\noindent As we can see from the example, the warping degree depends on the orientation. 
The pair $D_a$ is \textit{monotone} if $d(D_a)=0$. 
For example, $D_a$ depicted in Figure \ref{fig2} is monotone. 
Note that a monotone diagram is a trivial knot diagram.

\medskip
\hspace{3ex}The following theorem is our main theorem.

\phantom{x}
\begin{thm}
Let $D$ be an oriented knot diagram which has at least one crossing point. 
Then we have the following inequality: 
\begin{align*}
d(D)+d(-D)+1\leq c(D).
\end{align*}
Further, the equality holds if and only if $D$ is an alternating diagram. 
\label{mainthm}
\end{thm}
\phantom{x}

The rest of this paper is organized as follows. 
In section 2, we show some properties of the warping degree which are needed to prove Theorem \ref{mainthm}. 
In section 3, we prove Theorem \ref{mainthm} and show the table of warping degrees for all standard knot diagrams based on Rolfsen's table 
up to nine crossings \cite{rolfsen} and the warping degree of a standard torus knot diagram. 
In section 4, an application of Theorem \ref{mainthm} is given.

\section{Properties of the warping degree}

\medskip
\hspace{3ex} We show some interesting properties of the warping degree.
Let $-D$ be the inverse of $D$. 
Let $c(D)$ be the crossing number of $D$. 
Then we have the following lemma. 

\phantom{x}
\begin{lem}
Let $D_a$ be a diagram with a base point $a$. 
Then we have
\begin{align*}
d(D_a)+d(-D_a)=c(D).
\end{align*}
\end{lem}
\phantom{x}

\begin{figure}[h]
\begin{center}
\includegraphics[width=7em]{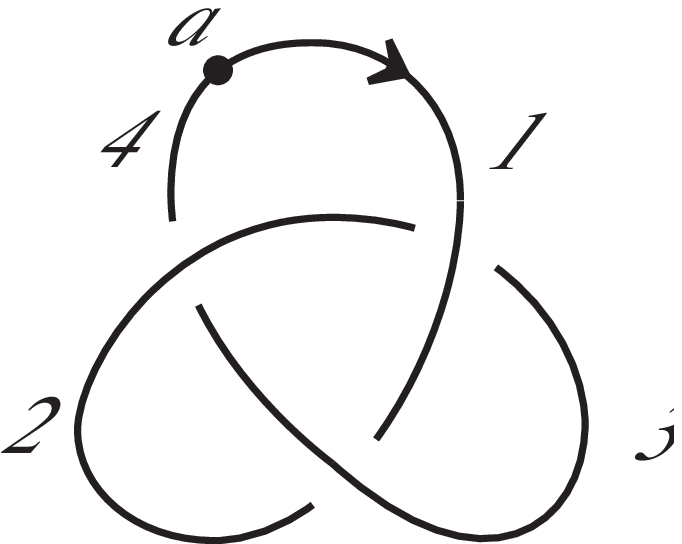}
\caption{}
\label{cut}
\end{center}
\end{figure}

Before proving Lemma 2.1, we introduce the method of judging locally whether a crossing point of an oriented knot diagram with a base point 
is a warping crossing point or not.
Let $D$ be an oriented knot diagram. 
Let $a$ be a base point of $D$. 
We notice that $D_a$ is divided into $c(D)+1$ arcs by cutting the base point and under-crossings. 
Then we label them along the orientation from the arc which has the base point as initial point as shown in Figure \ref{cut} .
Every crossing point consists of two or three arcs. 
We label each crossing point via the indices of the arcs in the following definition.

\phantom{x}
\begin{defn}
Let $p$ be a crossing point of $D_a$ which consists of an over-arc 
with the index $\alpha $ and the other two arcs with the index $\beta $, $\gamma $. 
We define the {\it cutting number} of $p$ in $D_a$, denoted by $cut_{D_a}(p)$, by the following formula:
\begin{align*}
cut_{D_a}(p)=2\alpha -\beta -\gamma .
\end{align*}
\end{defn}
\phantom{x}

\begin{figure}[h]
\begin{center}
\includegraphics[width=7em]{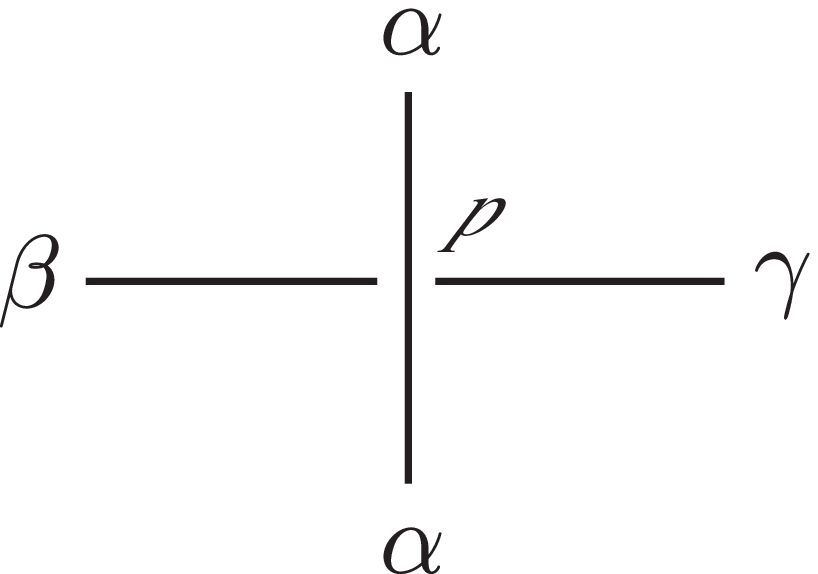}
\caption{}
\end{center}
\end{figure}

Note that the cutting number is always odd. 
Suppose $\beta $ is smaller than $\gamma $, namely $\gamma =\beta +1$. 
Then we have $cut_{D_a}(p)=2(\alpha -\beta )-1$ by substituting $\beta +1$ for $\gamma $.
Hence $cut_{D_a}(p)$ is odd. 
By the definition of the warping crossing point, we have the following lemma: 

\phantom{x}
\begin{lem}
A crossing point $p$ is a warping crossing point of $D_a$ if and only if $cut_{D_a}(p)>0$, 
and $p$ is a non-warping crossing point of $D_a$ if and only if $cut_{D_a}(p)<0$.
\end{lem}
\phantom{x}

We prove Lemma 2.1.

\phantom{x}

\noindent {\it Proof of Lemma 2.1.}
Let $(k_1, k_2, \dots , k_n)$ be the ordered set of arcs of $D_a$, where $k_i$ has the number $i$. 
Let $(l_1, l_2, \dots , l_n)$ be the ordered set of arcs of $-D_a$ as above.
Then we notice that $k_i$ and $l_{n+1-i}$ are the same arc except the orientations, for $i=1,2,\dots , n$.
Let $p$ be a crossing point of $D_a$ in Figure 4. 
Then we have the following equality:
\begin{align*}
cut_{-D_a}(p)&=2(n+1-\alpha )-(n+1-\beta )-(n+1-\gamma )\\
             &=-(2\alpha -\beta -\gamma )\\
             &=-cut_{D_a}(p).
\end{align*}

\noindent Therefore $p$ is a non-warping crossing point of $-D_a$ if and only if $p$ is a warping crossing point of $D_a$. 
Hence we obtain that
\begin{align*}
c(D)&=d(D_a)+\# \{ \text{non-warping crossing points of } D_a\} \\
    &=d(D_a)+\# \{ \text{warping crossing points of } -D_a\} \\
    &=d(D_a)+d(-D_a).
\end{align*}
This completes the proof. 
\hfill$\square$

\phantom{x}
We apply Lemma 2.1 to the mirror image in the following example. 

\phantom{x}
\begin{eg}
Let $D_a$ be an oriented knot diagram with a base point $a$, and $D_a^*$ the mirror image of $D_a$. 
Then we observe that a crossing point $p$ is a non-warping crossing point of $D_a^*$ if and only if 
$p$ is a warping crossing point of $D_a$. 
By Lemma 2.1, $p$ is a warping crossing point of $D_a^*$ if and only if $p$ is a warping crossing point of $-D_a$. 
Therefore we have $d(D^*)=d(-D)$. 
\end{eg}
\phantom{x}
  
For two base points which are put across a crossing point, we have the following lemma:

\phantom{x}
\begin{lem}
\begin{description}

\item[($1$)] For the base points $a_1$, $a_2$ which are put across an over-crossing in Figure 5, we have

\begin{align*}
d(D_{a_2})=d(D_{a_1})+1. 
\end{align*}

\begin{figure}[h]
\begin{center}
\includegraphics[width=7em]{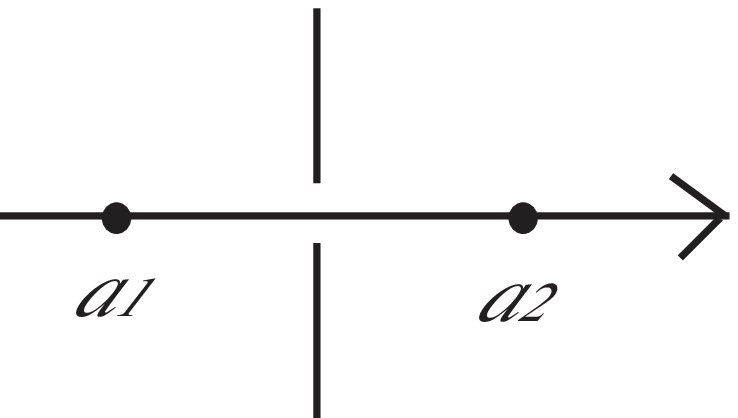}
\caption{}
\end{center}
\end{figure}
\phantom{x}

\item[($2$)] For the base points $a_1$, $a_2$ which are put across an under-crossing in Figure 6, we have

\begin{align*}
d(D_{a_2})=d(D_{a_1})-1. 
\end{align*}

\begin{figure}[h]
\begin{center}
\includegraphics[width=7em]{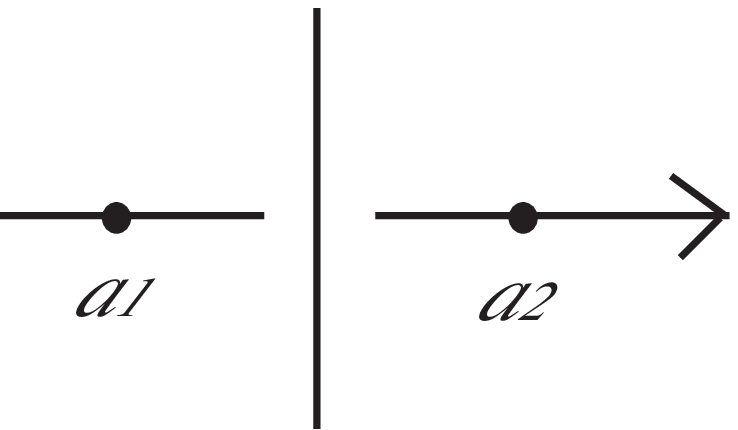}
\caption{}
\end{center}
\end{figure}

\end{description}
\end{lem}
\phantom{x}

By Lemma 2.5, we obtain the following lemma:

\phantom{x}
\begin{lem}
Let $D$ be an oriented alternating knot diagram.
Let $a$ be a base point of $D$ which is just before an over-crossing as shown in Figure 7, then we have the following equality:
\begin{align*}
d(D_a)=d(D).
\end{align*}

\begin{figure}[h]
\begin{center}
\includegraphics[width=8em]{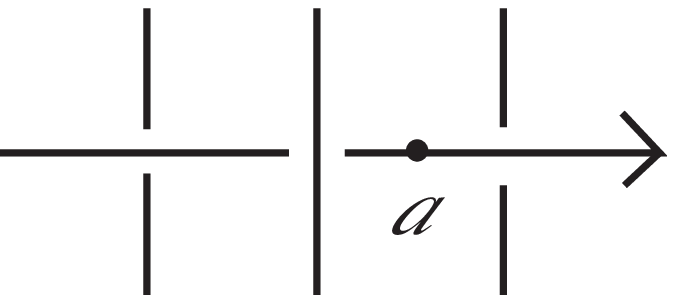}
\caption{}
\end{center}
\end{figure}
\end{lem}
\phantom{x}

\section{Proof of main theorem}

\phantom{x}

\hspace{3ex} We prove Theorem \ref{mainthm} in this section. 

\phantom{x}

\begin{figure}[h]
\begin{center}
\includegraphics[width=14em]{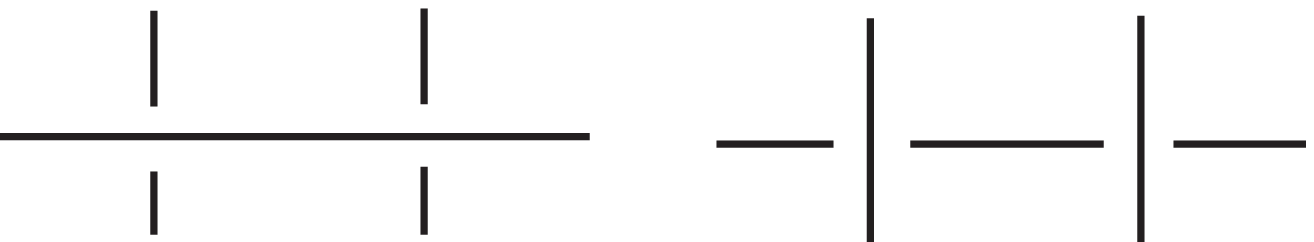}
\caption{}
\label{alt}
\end{center}
\end{figure}

\noindent {\it Proof of Theorem \ref{mainthm}.}
The following inequality holds for an arbitrary oriented knot diagram $D$ with $c(D)\geq 1$:
\begin{align}
\max _a d(D_a)-\min _a d(D_a) \geq 1.
\label{maxmin}
\end{align}
The equality holds if $D$ is an alternating diagram, by Lemma 2.5. 
On the other hand, if the equality holds, $D$ is an alternating diagram, 
namely there do not exist any two over-crossings or two under-crossings which are next to each other as shown in Figure \ref{alt} .

\noindent Let $a$ and $a'$ be base points which satisfy $d(D_a)=d(D)$ and $d(-D_{a'})=d(-D)$. 
Then we notice that 
\begin{align*}
\max _ad(D_a)=d(D_{a'}),
\end{align*} 
because $a'$ satisfies 
\begin{align*}
d(-D_{a'})=\min _a d(-D_a) 
\end{align*}
and a warping crossing point of $-D_a$ is a non-warping crossing point of $D_a$.
Hence the inequality (\ref{maxmin}) is equivalent to the following inequality:
\begin{align*}
d(D_{a'})-d(D_a) \geq 1,
\end{align*}
and this is equivalent to 
\begin{align*}
d(D_a)+1\leq d(D_{a'}).
\end{align*} 
By adding $d(-D_{a'})$ to each side, we have 
\begin{align*}
d(D_a)+d(-D_{a'})+1\leq d(D_{a'})+d(-D_{a'}). 
\end{align*}
By Lemma 2.1 and the conditions of $a$ and $a'$, we obtain the following inequality:
\begin{align*}
d(D)+d(-D)+1\leq c(D), 
\end{align*}
where the equality holds if and only if $D$ is an alternating diagram.
\hfill$\square$

\medskip
Here is an example of Theorem \ref{mainthm}.

\phantom{x}
\begin{eg}
This table lists all standard knot diagrams based on Rolfsen's knot table with crossing number 9 or less \cite{rolfsen}. 
In this table, $D(K)$ denotes the standard diagram of $K$ with the orientation which has the smaller warping degree, 
and a knot marked with \dag \ is non-alternating.
\begin{center}
\begin{tabular}{|p{5ex}|p{4em}|p{10ex}|p{42ex}|}

\hline
$K$ & $d(D(K))$ & $d(-D(K))$ & \\
\hline
$3_1$ & $1$ & $1$ & \\
\hline
$4_1$ & $1$ & $2$ & \\
\hline
$5_i$ & $2$ & $2$ & $i=1,2$ \\
\hline
$6_i$ & $2$ & $3$ & $i=1,2,3$ \\
\hline
$7_i$ & $3$ & $3$ & $i=1,2,\dots ,6$ \\
\hline
$7_7$ & $2$ & $4$ & \\
\hline
$8_i$ & $3$ & $4$ & $i=1,2,\dots ,17$ \\
\hline
$8_{18}$ & $2$ & $5$ & \\
\hline
$8_{19}$ \dag & $3$ & $3$ & \\
\hline
$8_{20}$ \dag & $2$ & $3$ & \\
\hline
$8_{21}$ \dag & $2$ & $2$ & \\
\hline
$9_i$ & $4$ & $4$ & $i=1,2,\dots ,13,16,18,20,21,23,25,27,28,$ \\
 & & & \hspace{3ex} $29,30,33,35,36,38,39,40$ \\
\hline
$9_i$ & $3$ & $5$ & $i=14,15,17,19,22,24,26,31,32,34,37,41$ \\
\hline
$9_i$ \hspace{1ex} \dag & $3$ & $3$ & $i=42,44,45,46$ \\
\hline
$9_i$ \hspace{1ex} \dag & $3$ & $4$ & $i=43,49$ \\
\hline
$9_{47}$ \dag & $2$ & $5$ & \\
\hline
$9_{48}$ \dag & $2$ & $3$ & \\
\hline
\end{tabular}\\
\end{center}
\end{eg}

\phantom{x}
In the following lemma, we determine the warping degree of the standard diagram of a torus knot. 

\phantom{x}
\begin{lem}
Let $T(p,q)$ be $(p,q)$-torus knot ($0<p<q$, $p$ and $q$ are coprime) and $D(p,q)$ the standard diagram of $T(p,q)$ 
with the orientation as shown in Figure \ref{fig9}. Then we have the following: 

\begin{description}
\item[(1)] $d(D(p,q))=d(-D(p,q))=\frac{(p-1)(q-1)}{2},$ 
\item[(2)] $c(D(p,q))-d(D(p,q))-d(-D(p,q))=p-1.$ 
\end{description}

\begin{figure}[h]
\begin{center}
\includegraphics[width=12em]{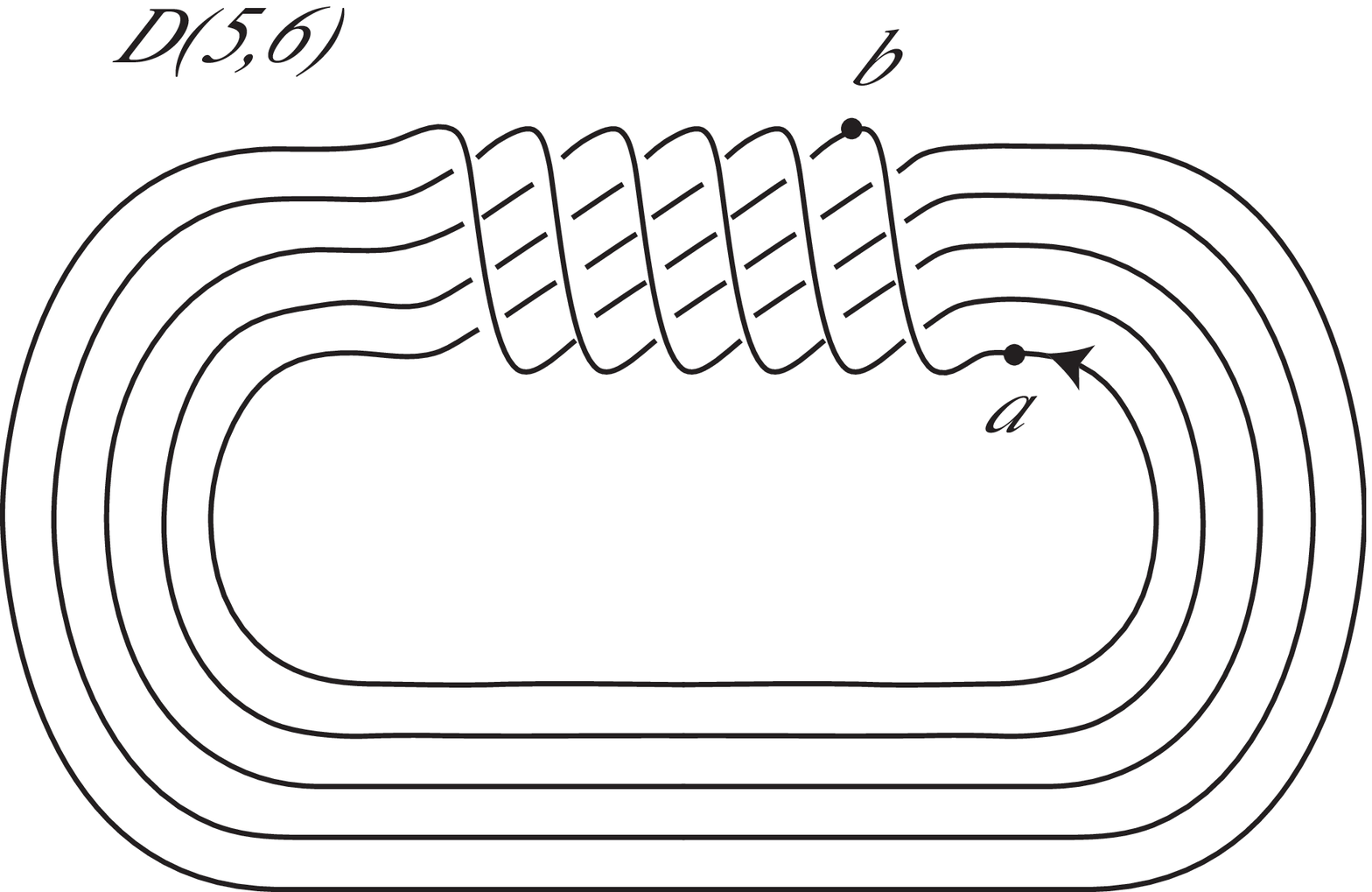}
\caption{}
\label{fig9}
\end{center}
\end{figure}

\begin{proof}
In $D(p,q)$, there are $(p-1)$ over-crossings and $(p-1)$ under-crossings alternately. 
By Lemma 2.5, a base point which is just before $(p-1)$ over-crossings realizes the warping degree of the diagram. 
For example, base points $a$ and $b$ in Figure \ref{fig9} satisfy 
\begin{align*}
d(D(p,q)_a)=d(D(p,q)),\\
d(-D(p,q)_b)=d(-D(p,q)),
\end{align*}
respectively. 
Considering the upside-down image of $D(p,q)$, we have 
\begin{align*}
d(D(p,q))=d(-D(p,q)), 
\end{align*}
that is, 
\begin{align}
d(D(p,q)_a)=d(-D(p,q)_b). 
\end{align}
With respect to the crossing number of the diagram, we have 
\begin{align}
d(D(p,q)_a)+d(-D(p,q)_a)=c(D)=(p-1)q 
\end{align}
by Lemma 2.1. 
And by Lemma 2.5, we have the following relation 
\begin{align}
d(-D(p,q)_b)=d(-D(p,q)_a)-(p-1). 
\end{align}
By the formula (2) and (4), we have 
\begin{align}
d(D(p,q)_a)=d(-D(p,q)_a)-(p-1). 
\end{align}
And by the formula (3) and (5), we have 
\begin{align*}
2d(D(p,q)_a)=(p-1)q-(p-1). 
\end{align*}
Hence we obtain the warping degree
\begin{align*}
d(D(p,q))=\frac{(p-1)(q-1)}{2}. 
\end{align*}
Hence we have 
\begin{align*}
c(D(p,q))-d(D(p,q))-d(-D(p,q))\\
=(p-1)q-\frac{(p-1)(q-1)}{2}\times 2=p-1.
\end{align*}
This means that $c(D)-d(D)-d(-D)$ depends only on $p$ in this case. 
\end{proof}
\end{lem}\par
\phantom{x}

\section{To a knot invariant}

We apply Theorem \ref{mainthm} to knots. 
For a knot $K$ we define an invariant $e(K)$ of $K$ by the following formula:
\begin{align*}
e(K)=\min \{ d(D)+d(-D) | D: \text{an oriented diagram of} \ K \text{ with } c(D)=c(K) \} .
\end{align*}
Then we obtain the following theorem: 

\phantom{x}
\begin{thm}
Let $K$ be a non-trivial knot, and $c(K)$ the crossing number of $K$. Then we have the following (1) and (2). 
\begin{description}
\item[(1)] We have an inequality $e(K)+1\leq c(K).$ 
Further, the equality holds if and only if $K$ is a prime alternating knot.
\item[(2)] For any positive integer $n$, there exists a prime knot $K$ such that $$c(K)-e(K)=n.$$ 
\end{description}

\medskip
\begin{proof}
First, we show the equality of (1). 
By Theorem \ref{mainthm}, we have the equality 
\begin{align}
d(D)+d(-D)+1=c(D)
\label{alt31}
\end{align}
if and only if $D$ is an alternating diagram. 
If $K$ is a prime alternating knot, 
then minimal crossing diagrams of $K$ are alternating \cite{murasugi}. 
Hence we have the equality by considering the minimum of the equality (\ref{alt31}). 
On the other hand, if $K$ is a non-prime alternating knot, 
then there is a minimal crossing non-alternating diagram so that $e(K)+1<c(K)$ \cite{menasco}. 
Secondly, we look to $c(T(p,q))-e(T(p,q))$ for the $(p,q)$-torus knot $T(p,q)$
$(0<p<q)$ to prove (2). 
Schubert mentioned in \cite{schubert} (cf.\cite{murasugi2}) that 
\begin{align*}
c(T(p,q))=(p-1)q. 
\end{align*}
Ozawa showed in \cite{ozawa} that the ascending number of $T(p,q)$, 
which is equal to the minimal warping degree for all diagrams of $T(p,q)$ and all orientations, 
is $(p-1)(q-1)/2$. 
Then we have 
\begin{align*}
e(T(p,q))=\frac{(p-1)(q-1)}{2}+\frac{(p-1)(q-1)}{2}=(p-1)(q-1), 
\end{align*}
because $d(D(p,q))=d(-D(p,q))$. Hence we have 
\begin{align*}
c(T(p,q))-e(T(p,q))=(p-1)q-(p-1)(q-1)=p-1.
\end{align*}
\end{proof}
\end{thm}
\phantom{x}

\section*{Acknowledgments}
The author would like to thank Akio Kawauchi and Taizo Kanenobu for their kindness in teaching her. 
She also would like to thank Oyuki Hermosillo, Ryo Hanaki, Ryo Nikkuni, 
and the members of Friday Seminar on Knot Theory 
in Osaka City University for valuable conversations, encouragement, and feedback on previous version of the paper.

\maketitle


\begin{thebibliography}{9}
\bibitem{kawauchi}A.~Kawauchi, Lectures on knot theory (in Japanese), Kyoritsu Shuppan (2007).
\bibitem{menasco}W.~Menasco, Closed incompressible surfaces in alternating knot and link complements, Topology 23 (1984), 37--44. 
\bibitem{murasugi}K.~Murasugi, Jones polynomials and classical conjectures in knot theory, Topology 26 (1987), 187--194.
\bibitem{murasugi2}K.~Murasugi, On the braid index of alternating links, Trans. Amer. Math. Soc. 326 (1991), 237--260. 
\bibitem{ozawa}M.~Ozawa, Ascending number of knots and links, arXiv:math.~GT/07053337v1.~(2007).
\bibitem{rolfsen}D.~Rolfsen, Knots and links, Publish or Perish, Inc. (1976).
\bibitem{schubert}H. Schubert, {\"U}ber einenumerische Knoteninvariante, Math.~Zeit.~61(1954), 245--288.
\end{thebibliography}
\end{document}